\documentstyle[11pt]{article}
\begin{document}\date{}
\title{Prime rings with PI rings of constants}
\author{V. K. Kharchenko \thanks{The first author wishes to thank 
CONACYT--M\'exico for its support, Catedra Patrimonial number 940411-R96 and also 
Russian Found of Fundamental Research, grant 95-01-01356.}
\\ J.Keller  \\ S. Rodr\'\i guez-Romo \thanks{ The third author wishes to thank
CONACYT--M\'exico for its support under grant 4336E9406.}
\\ Centre of Theoretical Research
\\ UNAM, Campus Cuautitl\'an
\\ Apdo. Postal 95, Cuautitl\'an Izcalli
\\ Estado de M\'exico, 54768, M\'exico }
\maketitle

{\it in memory of S.A.Amitsur}

\begin{abstract}It is shown that if the ring of constants of a restricted 
differential Lie algebra with a quasi-Frobenius inner part satisfies a
polynomial identity (PI) then the original prime ring has a generalized
polynomial identity (GPI). If additionally the ring of constants is 
semiprime then the original ring is PI. The case of a
non-quasi-Frobenius inner part  is also considered.
\end{abstract}

\section{Introduction}

Rings of constants of restricted differential Lie algebras with an outer
action on prime and semiprime rings were investigated in detail in 
papers [Kh82], [Po83], [Pi86] (see also [Kh91,Ch.4, Ch6(6.4)]). In the 
present paper we are going to consider actions with a nontrivial inner
part. In the papers [Ko91] and [Kh81] it is shown that the minimal restriction
required is that the inner part should be quasi-Frobenius (selfinjective).
We are interested in the structure of a prime ring $R$ provided it is known 
that its ring of constants satisfies a polynomial identity. I.V.L'vov's
example [Lv93] shows that in this case the ring $R$ does not need to
be a  PI-ring. We will show that in this case $R$ satisfies a generalized
polynomial identity.

The notion of a generalized polynomial identity was introduced by 
S.A.Amitsur in [Am65]. In his paper S.A.Amitsur proved a structure 
theorem for primitive rings with generalized polynomial identities.
Later W.S.Martindale [Ma69] generalized this result to arbitrary
prime rings. Using this theorem 
we will prove that if the ring of constants is a semiprime PI-ring
and the inner part is quasi-Frobenius, then the ring $R$ is a PI-ring.

\section{Preliminaries}

Recall that a {\it derivation} of a ring $R$ is an additive mapping
$d:R \rightarrow R$ satisfying the condition
$(xy)^d=x^dy+xy^d.$
If $d_1, d_2$ are derivations then it is easy to see that the
commutator $[d_1,d_2]=d_1d_2-d_2d_1$ is also a derivation. Therefore
the set $Der R$ of all derivations of $R$ is a Lie subring in the ring
of endomorphisms of the abelian group $(R,+)$. Moreover, if $z$ is a central
element, then the composition of $d$ with the multiplication by $z$
is a derivation
$$(xy)^{dz}=z(xy)^d=(zx^d)y+x(zy^d)$$
In this case the operators of multiplication may not commute with derivations:
$x^{zd}\stackrel{\rm def}{=} (zx)^d=z^dx+zx^d$ or
$$zd=dz+z^d. \eqno(1)$$

Thus the set $Der R$ is a right module over the center $Z.$ The 
module structure of $Der R$ 
is connected with the commutator operation by the formula
$$[dz,d_1]=[d,d_1]z+dz^{d_1}. \eqno(2)$$
Note that $z^{d_1}$ is again a central element:
$[z^{d_1},x]=[z^{d_1},x]+[z,x^{d_1}]=[z,x]^{d_1}$=0.

Finally, if the characteristic $p$ of the ring $R$ is nonzero, $pR=0,$
then the $p$th power of any derivation will be a derivation by the Leibniz
formula
$$(xy)^{d^p}=\sum_{k=0}^{k=p} C_p^k x^{d^k}y^{d^{p-k}}=x^{d^p}y+xy^{d^p}.$$
Now it is natural to formulate the following definition.

{\bf 2.1. Definition.} A set of derivations is called a {\it differential
restricted Lie $Z$-algebra}, or shortly a {\it Lie $\partial $-algebra}, if
it is a right $Z$-submodule of $DerR$ closed with respect to the operations
 $[d_1,d_2]=d_1d_2-d_2d_1$ and $d^{[p]}=d^p.$

 Note that the notion of a Lie $\partial $-algebra can be formalized
  abstractly as a  restricted Lie ring with a structure of right
 $Z$-module connected with the main operations by formula (2) and the
following formula
 $$(dz)^{[p]}=d^{[p]}z^p+d\cdot
 (\ldots ((\overbrace{z^dz)^dz)^d\ldots )^d}^{p-1}z  \eqno(3)$$
 which follows from (1) (see details in [Kh91, pp. 6-11]; for a slightly more
  general approach see in [Pa87]).

 Now let $R$ be a prime ring. Denote by $R_{\cal F}$ its left Martindale ring
of quotients (see, for example, [Kh91 pp.19-24]), by $Q$ the symmetric Martindale
ring of quotients. Recall that the center $C$ of $R_{\cal F}$ is called the
{\it extended} (or {\it generalized}) centroid of $R$ and it is a field (see [Ma69]).
All derivations of $R$ can be uniquely extended to derivations of $Q$ and
 of $R_{\cal F}.$ The extended derivations are characterized in $Der Q$ by the
property $R^d\subseteq R$ but  the linear combinations over $C$ of extended
derivations do not satisfy this property. Therefore we have to consider more
general objects.

{\bf 2.2. Definition.} A derivation $d$ of $Q$ is called $R$-{\it continuous}
if there
exists a nonzero two-sided ideal $I$ of $R$ such that $I^d\subseteq R.$

It is easy to see that the set ${\cal D}(R)$ of all $R$-continuous derivations
 is a differential restricted Lie $C$-subalgebra of $Der Q.$

In the present paper we consider Lie
$\partial $-algebras  of $R$-continuous derivations which are finite
dimensional over $C.$

Let us fix the notations $R, C, Q, R_{\cal F}, {\cal D} (R)$
for a prime ring, its extended centroid, the symmetric Martindale ring of
quotients, the left Martindale ring of quotients
 and the Lie $\partial $-algebra of $R$-continuous derivations,
respectively. Throughout the paper $L$ denotes a
restricted differential Lie $C$-algebra of $R$-continuous derivations,
$L\subseteq {\cal D} (R),$ finite dimensional over $C,$
 and $R^L=\{ r\in R:\forall \mu \in L \ \
\ r^{\mu }=0\} $ is its ring of constants.

\section{The inner part of a Lie $\partial $-algebra}

If $a$ is an element of $Q$ then the map $a^{-}:x\rightarrow xa-ax$
is an $R$-continuous derivation, i.e. $Q^{-}\subseteq {\cal D}(R).$

{\bf 3.1. Definition.} The space $K(L)$ generated over $C$ by all
$q\in Q$ such that $q^{-}\in L$ is called the {\it inner linear part of $L.$}

It is clear that $C^{-}=0,$ therefore $K(L)$ contains $C$ and in particular
it contains the unit of $Q.$

{\bf 3.2. Lemma.} {\it The space $K(L)$ is a restricted Lie subalgebra of the
adjoint restricted Lie algebra $Q^{(-)}.$}

Recall that $Q^{(-)}$ is a restricted Lie algebra defined on the $C$-space
$Q$ with the operations $[q_1,q_2]=q_1q_2-q_2q_1,\ \  q^{[p]}=q^p.$

For the proof of the lemma it is enough to show that $K(L)$ is closed
with respect to these operations. This fact immediately follows from the
formulae
$$[a,b]^{-}=[a^{-},b^{-}] \eqno(4)$$
$$(a^p)^{-}=(a^{-})^{[p]}. \eqno(5)$$

{\bf 3.3. Lemma.} {\it $K(L)^{-}$ is equal to the subalgebra $L_{int}$
of all inner derivations of $L.$}

The proof is evident.

{\bf 3.4. Definition.} The associative subalgebra ${\cal B}(L)$ 
 generated in $Q$ by $K(L)$ is called the {\it inner associative part of} $L.$

{\bf 3.5. Lemma.} {\it The algebra ${\cal B}(L)$ is of finite dimension
over $C.$}

{\bf Proof.} By the definition of operations in $K(L),$ the identity map
$id$ is a homomorphism of restricted Lie algebras
$id: K(L)\rightarrow {\cal B}(L)^{(-)}.$ Therefore ${\cal B}(L)$ as an associative
envelope  of ${\cal B}(L)^{(-)}$ is a homomorphic image of the universal
restricted associative envelope $U_p(K(L)).$ The latter has dimension
$(\dim K(L))^p.$ The lemma is proved.

{\bf 3.6. Lemma.} {\it The algebra ${\cal B}(L)$ is stable under the action of
$L,$ i.e. ${\cal B}(L)^{\mu } \subseteq {\cal B}(L)$ for all} $\mu \in L.$

The proof  follows from the formula
$$(q^{\mu })^{-}=[q^{-},\mu ]. \eqno(6)$$

\section{Differential operators}

Denote by $\Phi (L)$ the associative 
subring generated in the endomorphism ring
of the abelian group $(Q,+)$ by $L$ and by the operators of left and
right multiplications by elements from ${\cal B}(L).$ By formula
(1) the ring $\Phi (L)$ may not be an algebra over $C.$ Of course
$\Phi (L)$ is an algebra over the subfield of central constants
$$F=C^L\stackrel{\rm def}{=} \{ c\in C:\forall l\in L \ \  c^l=0 \} .$$
Nevertheless $\Phi (L)$ is a left and a right space over $C$ 
while the subring of left multiplications, ${\cal B}(L)^l,$ and that of
right multiplications, ${\cal B}(L)^r,$ are algebras over $C.$

{\bf 4.1.} Let us fix derivations $\mu _1, \ldots ,\mu _m \in L$ such
that $\mu _1+K(L)^{-}, \ldots ,\mu _m+K(L)^{-}$
form a basis for the right $C$-space $L/K(L)^{-}.$ An operator $\Delta $
is called {\it correct} if it is of the form:
$$\Delta =\mu ^{s_1}_1\mu ^{s_2}_2 \ldots \mu ^{s_m}_m,$$
where $0\leq s_i<p$ and we suppose that $\mu ^0=1$ is the identity
operator.

Let $U$ be a right linear space generated by all correct operators. By
formula (1) this set will be a left space over $C,$ also.

{\bf 4.2. Proposition.} {\it The ring $\Phi (L)$ of differential operators
 is isomorphic as a left and a right space over $C$ to a
tensor product over $C:$
$$\Phi (L)\simeq {\cal B}(L)^r\otimes {\cal B}(L)^l\otimes U\simeq
U\otimes {\cal B}(L)^l\otimes {\cal B}(L)^r,  \eqno(7)$$
where $U$ is the linear space generated by correct operators over $C$.}

{\bf Proof.} It is enough to show that each differential operator
$d\in \Phi (L)$ has a unique representation in the form
$$d=\sum _{i,j,k} \alpha ^{(k)}_{ij}a^r_{ik}a^l_{jk}\Delta _k  \eqno(8)$$
and a unique representation in the form
$$d=\sum _{i,j,k} \Delta _ka^l_{ik}a^r_{jk}\alpha ^{(k)}_{ij},  \eqno(9)$$
where $a_{ik}, a_{jk} \in A$ and $A$ is some fixed basis 
of ${\cal B}(L)$ over $C$
(recall that by associativity, $a^r_{ik}a^l_{jk}=a^l_{jk}a^r_{ik}$)
and the $\Delta _k$'s
are correct words in $\{ \mu _1,\ldots ,\mu _m \} .$

The existence of this presentation follows from the relations
$$\mu a^r=a^r\mu -(a^{\mu })^r  \eqno(10)$$
$$\mu a^l=a^l\mu -(a^{\mu })^l  \eqno(11)$$
$$\mu ^p=\mu _1c_1+ \ldots +\mu _mc_m+b^r-b^l \eqno(12)$$
$$\mu _i\mu _j=\mu _j\mu _i+\mu _1c_1+ \ldots +\mu_mc_m +b^r-b^l,  \eqno(13)$$
where in formula (12) \ $\mu _1c_1+ \ldots +\mu _mc_m+b^{-}$ \  
is a representation
of $\mu ^p\in L$ as a linear combination of $\mu _i$'s modulo $K(L)^{-}$
and in (13) \  $\mu _1c_1+ \ldots +\mu _mc_m+b^{-}$ \  is the
corresponding representation
of $[\mu _i,\mu _j]\in L.$

The transformations of the left hand sides to the right hand sides
 (in the last formula
only if $i>j$) allow us to reduce the operator to the form (8).

If we write formulae (10), (11) in the form
$$a^r\mu =\mu a^r+(a^{\mu })^r \eqno(14)$$
$$a^l\mu = \mu a^l+(a^{\mu })^l \eqno(15)$$
then in the same way the operator is reduced to the form (9).

For the proof of the uniqueness it is possible to use
the following results on differential
identities (see [Kh91, theorem 2.2.2, corollary 2.5.8] or [Kh78]).

{\bf 4.3. Proposition.} {\it If the derivations $\mu _1, \ldots ,\mu _m \in
{\cal D}(R)$ are linearly independent modulo $Q^{-},$ and if
the ring $R$ satisfies an identity of the type
$$\sum ^{p^n}_{k=1} \sum _i a_{ki} x^{\Delta _k}b_{ki}=0,$$
where $\Delta _1, \ldots , \Delta _{p^n}$ --- are all correct operators and the
coefficients $a_{ki},b_{ki}$ belong to $R_{\cal F},$ then
$\sum _ia_{ki} \otimes b_{ki}=0$ in $R_{\cal F} \otimes _C R_{\cal F}$
for all $k, 1\leq k \leq p^n.$ In the same way if the identity
$$\sum ^{p^n}_{k=1} (\sum _i a_{ki}xb_{ki})^{\Delta _k}=0$$
is valid then $\sum _i a_{ki}\otimes b_{ki}=0, 1\leq k\leq p^n.$}

Since ${\cal D}(I)={\cal D}(R)$ and $Q(I)=Q(R)$ for each nonzero ideal $I$
of $R$ (see [Kh91, Lemma 1.8.4]), then proposition 4.3 shows that the restriction of a nonzero
differential operator $d\in \Phi (L)$ to $I$ is nonzero. This note is
important due to the following lemma:

{\bf 4.4. Lemma.} {\it For each differential operator $d\in \Phi (L)$
there exists a nonzero two sided ideal $I$ of $R$ such that
$I^d\subseteq R.$}

The proof is easily  obtained by induction from the formula
$(I^2)^{\mu }=I^{\mu }I+II^{\mu }\subseteq I$
which is valid for the ideal $I$ such that $I^{\mu }\subseteq R.$

\section{Quasi-Frobenius algebras}

Recall that a finite dimensional algebra $B$
over a field $C$ is called quasi-Frobenius if
one of the following equivalent conditions is valid (see [CR62]).

(Q1) {\it For each left ideal $\lambda $ and right ideal $\rho $ of $B$
the following equalities hold:
$$l(r(\lambda ))=\lambda ,\ \ \  r(l(\rho ))=\rho ,$$
where $l(A)=\{ b\in B: bA=0\} $ is the left annihilator,
$r(A)=\{ b\in B: Ab=0\} $ is the right annihilator.}

(Q2){\it  The left regular module $_BB$ is injective.}

(Q3) {\it Modules $_BB$ and $(B_B)^*=Hom (B,C)$ have the same indecomposable
components.}

Recall that for any left (right) module $M$ the set of  all linear
functionals $M^*$ has a structure of right (left) module defined by the
formula $(m^*b)(m)=m^*(bm)$ (respectively $m(bm^*)=(mb)m^*$).
The modules $M$ and $N$ for $N\simeq M^*$ are called 
{\it conjugated modules}. If the
module $M$ is of finite dimension then $(M^*)^*\simeq M$ and the
conjugacy of modules (left and right), $M$ and $N,$ can be
characterized by the existence of a nondegenerate  associative bilinear form
$(\ ,\ ):N\times M\rightarrow C.$ In this case for every basis
$a_1,\ldots ,a_n$ of $M$ there exists a {\it dual} basis 
$a_1^*,\ldots ,a_n^*$ of $N$ which is characterized by the
following properties
$(a_i^*,a_i)=1, (a_i^*,a_j)=0, i\neq j.$ 

Condition (Q3) implies the following condition which is important for us:

(Q4) {\it The sum of all right ideals $\rho $ of $B$
conjugated to left ideals of $B$ is equal to
$B.$}

It can be proved that this condition is also equivalent to $B$ being
quasi-Frobenius. Moreover, as (Q1) is left-right symmetric then the left
analog of (Q4) is also valid.

(Q5) {\it The sum of all left ideals ${\lambda }$ of $B$ conjugated to right
ideals of $B$ is equal to $B.$}

The most important subclass of the class of quasi-Frobenius algebras is
the class of Frobenius algebras. These algebras are defined by one of the
following equivalent conditions ([CR62]).

(F1) {\it For each left ideal $\lambda $ and right ideal $\rho $ of $B$
 the following equalities hold:
 $$l(r(\lambda ))=\lambda ,\ \ \dim r(\lambda )+\dim \lambda =\dim B$$
\ $$ r(l(\rho ))=\rho , \ \ \dim l(\rho )+\dim \rho=\dim B.$$ }

(F2) {\it There exists an element 
$\varepsilon \in B^*$ whose kernel contains
no nonzero onesided ideals of $B.$}

(F3) {\it There exists a nondegenerate associative bilinear form
$B\times B\rightarrow C.$  }

(F4) {\it The modules $_BB$ and $(B_B)^*$ are isomorphic.}

Classical examples of Frobenius algebras are: group algebras of finite
groups over a field of arbitrary characteristic, universal restricted 
enveloping algebras of finite dimensional Lie $p$-algebras, finite dimensional 
Hopf algebras, Clifford algebras. Finite dimensional semisimple algebras 
evidently satisfy (F1) therefore they are Frobenius.

\section{Universal constants}

Let $\lambda $ and $\rho $ be left and right conjugated ideals of
${\cal B}(L).$ Let us choose a basis $a_1, \ldots ,a_n$ of $\lambda $
and let $a_1^*, \ldots ,a_n^*$ be the dual basis of $\rho .$
It is well-known that the element $c=\sum a_i\otimes a_i^*$ of  the
tensor product ${\cal B}\otimes _C{\cal B}$ commutes with the elements of
${\cal B}, \ \ bc=cb$ for all $b\in {\cal B}.$ This implies that the
set of  values of the operator $c_{\lambda ,\rho }=\sum a_i^l(a_i^*)^r$
is contained in the centralizer of ${\cal B}.$ In particular for any
$\mu \in K(L)^-$ we have $$c_{\lambda ,\rho }(x)^{\mu }=0. \eqno(16) $$

Let $U(L)$ be the associative subring of $\Phi (L)$ generated by $L$
and by the operators of multiplication by central elements. It is clear that
$U(L)$ is both a left and a right space over $C$ and an algebra over the field
of central constants $F=C^L.$

Consider the right ideal $I=K(L)^-\cdot U(L)$ of $U(L).$ 
First of all the formula
$\mu a^-=a^-\mu -(a^{\mu })^-$ shows that $I$ is a two sided ideal of $U(L).$
The same formula and  formulae (12), (13)
show that the identity operator and operators of the
form $a_1^-a_2^-\ldots a_s^-\Delta ,$
where $\Delta $ is a correct operator, $a_i\in K(L), \ s\geq 0,$
generate $U(L)$ as a left space over $C.$

{\bf 6.1. Proposition.} {\it The factor-algebra $U(L)/I=\overline{U}$
 is Frobenius as
an algebra over $F=C^L.$}

{\bf Proof.} By the well-known R.Baer theorem [Ba27] the dimension of $C$
over $F$ is finite and therefore $\overline{U} =U(L)/I$ has a finite dimension
over $F.$ Since $K(L)^-\subseteq I,$ the elements $\bar{\mu }_1=\mu _1+I,
\ldots , \bar{\mu }_m=\mu _m+I$ generate $\overline{U} $ as a ring over $C.$
Moreover the relations $\bar{\mu } _i\bar{\mu }
_j=\overline{[\mu _i,\mu _j]}+\bar{\mu } _j\bar{\mu } _i$ show that the
 images of
correct words
$\bar{\Delta } _k$ generate $\overline{U}$ as a left vector space over $C.$
The main note is that the elements $\bar{\Delta } _k$
 are linearly independent over $C.$ If
 $$\sum _kc_k\Delta _k=\sum _k d_k\Delta _k\in I,$$
 where $d_k$ are linear combinations of products of the type
 $a_1^-\cdots a_s^-,$ then taking into account that $a^-=a^r-a^l$
 and using Proposition 4.3, we have $c^r_k=d_k$ for all $k,$
 which is impossible since $c_k^r(1)=c_k,\ d_k(1)=0.$
 Thus $\bar{\Delta }_k$ are linearly independent.

Now let us define Berkson's linear map
(see [Be64]) $\varphi :\overline{U} \rightarrow C$ which
 corresponds to
the element $\sum c_k\bar{\Delta }_k$
the coefficient of $\bar{\Delta }_{p^m}=\bar{\mu }^{p-1}_1\ldots \bar{\mu
}^{p-1}_m.$ The kernel of this linear map contains neither left nor right
nonzero ideals, since the product
$$(\bar{\mu }_1^{s_1}\ldots \bar{\mu }_m^{s_m})(\bar{\mu }^{p-s_1-1}_1\ldots
\bar{\mu }^{p-s_m-1}_m)$$ written as a linear combination of correct words
contains a unique member $\bar{\Delta }_{p^m}$ with a coefficient
equal to 1.

If $\psi :C\rightarrow F$ is any projection, then the linear functional
$\varepsilon :d\mapsto \psi (\varphi (d))$ satisfies (F2) and therefore
$\overline{U}$ is a Frobenius algebra. The proposition is proved.

Let us consider the right subspace $\hat{U}$ of $\overline{U}$ over $C$
generated by all nonempty words $\bar{\Delta }_k.$ This space does not
contain the unit (the identity operator) and it is a right 
(but, possibly, not
a left) ideal because by  formula (14) one has
$$\bar{\Delta }c\bar{\mu }=\bar{\Delta }\bar{\mu }c+\bar{\Delta }c^{\mu }.$$
By formula (13), the product $\bar{\Delta }\bar{\mu }$  can be
written as a linear combination $\sum \bar{\Delta }_kc_k,$ where
$\bar{\Delta }_k$ are nonidentity correct operators.

Thus, the left annihilator $A=l(\hat{U})$ in the algebra $\overline{U}$ is not
equal to zero. Moreover, by (F1) its dimension over $F$ is connected with
 the dimension of $\hat{U}$ by the formula
 $ \dim _F\overline{U}=\dim _F\hat{U}+\dim _FA.$ On the other hand
 $\dim _F\overline{U}=\dim _F\hat{U}+\dim _FC$ i.e. the dimensions of $A$
 and $C$ over $F$ coincide. It means that $A$ is one dimensional over $C$
 i.e. $A=C\bar{f}$ (but possibly $A\neq \bar{f}C$ as $A$ may  not be a
 right $C$-space), where $\bar{f}=\sum \bar{\Delta }_kc_k=\sum
 c_k^{\prime }\bar{\Delta }_k$ is a nonzero element of $\overline{U}.$

Thus, we have obtained that $\bar{f}\bar{\mu }_i=\bar{0}$ in $\overline{U}.$
In the ring of differential operators this means that
$f\mu _i\in K(L)^{-}\cdot U(L).$ We have also that
$fK(L)^{-}\subseteq K(L)^{-}\cdot U(L)$ as $I=K(L)^{-}\cdot U(L)$
is a two sided ideal. Thus
$$fL\subseteq f(\sum (\mu _iC+K(L)^{-}))\subseteq K(L)^{-}\cdot U(L)$$
which, using formula (16), implies
$$((c_{\lambda ,\rho } (x))^f)^{\mu }=0 \eqno(17)$$
for all ${\mu \in L}.$ Let us formulate the obtained result as a lemma
(see also Lemma 4.6, [Kh95]).

{\bf 6.2. Lemma.}  {\it There exists a differential operator $f$ of the type
$\sum \Delta _kc_k=\sum c^{\prime }_k\Delta _k,$ such that for each
conjugated left ideal $\lambda $ and right ideal $\rho $ of ${\cal B}$
with dual bases $a_1, \ldots , a_n$ and $a_1^*,\ldots ,a_n^*,$ the
operator
$$u_{\lambda , \rho }=\sum_i a_i^l(a_i^*)^rf \eqno(18)$$
has values only in the ring of constants $Q^L.$ There exists a nonzero
ideal $I$ of $R$ such that
$$0\neq u_{\lambda ,\rho }(I)\subseteq R^L \eqno(19)$$}

{\bf Proof.} The representation of 
$f$ in the form $\sum c^{\prime }_k\Delta _k$
follows from (10).
Formula (19) follows from formula (17), proposition 4.3 and lemma 4.4.

\section{PI rings of constants}

In this secton we will prove the theorem about a generalized polynomial
identity and discuss
its generalization to the case when the inner part is not quasi-Frobenius.

{\bf 7.1. Theorem.} {\it Let $L$ be a finite dimensional restricted
differential Lie $C$-algebra of $R$-continuous derivations of a prime
ring $R$ of positive characteristic $p>0.$ Suppose that the inner associative
part ${\cal B}(L)$ of $L$ is quasi-Frobenius. If the ring of constants $R^L$
is PI then $R$ is GPI.}

{\bf Proof.} Let $f(x_1,\ldots ,x_n)=0$ be a multilinear identity
of $R^L.$ Let us choose arbitrary left ideals $\lambda _1,\ldots ,\lambda _n$
of ${\cal B}(L)$ having conjugated right ones $\rho _1, \ldots ,\rho _n.$
By Lemma 6.2 for every $j, 1\leq j\leq n$ there exists an operator
$$u_j=u_{\lambda _j,\rho _j}=\sum_i a^l_{ij}(a_{ij}^*)^rf_j=\sum_{i,k}
a_{ij}^l(a_{ij}^*)^rc_k^{\prime }\Delta _k$$
and a nonzero ideal $I_j$ of $R,$ such that $0\neq u_j(I_j)\subseteq R^L.$
If $I=\cap I_j$ then $u_j(x)\in R^L$ for all $x\in I$
and therefore the following differential identity holds in $I$ \
$$f(u_1(x_1),\ u_2(x_2),\ldots ,\ u_n(x_n))=0.$$

Let us fix some values of $x_2=b_2, \ldots ,x_n=b_n$ in $I.$ We have
$$f(\sum_{i,k} (c_k^{\prime }a_{i1}x_1a_{i1}^*)^{\Delta _k},\  u_2(b_2),\ldots
,\ u_n(b_n))=0. \eqno(20)$$
By Leibnitz formula any expression of the type $(axb)^{\Delta }$
 can be written in the form
 $$(axb)^{\Delta }=ax^{\Delta }b+\sum_s a_sx^{\Delta _s}b_s,$$
 where $\Delta _s$ are subwords of $\Delta .$ In particular
 $$(c_k^{\prime }a_{i1}x_1a_{i1}^*)^{\Delta _k}=
 c_k^{\prime }a_{i1}x_1^{\Delta _k}a_{i1}^*+\sum_s a_sx_1^{\Delta _s}b_s.
  \eqno(21)$$
  If $\Delta _{k_0}$ is the greatest operator such that $c^{\prime }_{k_0}$
  is not zero, then this formula allows us to represent (20) in the form
  $$\sum^{k_0}_{k=1} \sum_i v_{ki}x_1^{\Delta _k}w_{ki}=0,$$
  here we suppose that $\Delta _1<\Delta _2<\ldots <\Delta _{p^m}$ is the
  lexicographic ordering of all correct operators.
  By Proposition 4.3 applied to the prime ring $I$ we have
  $$\sum_i v_{k_0i}\otimes w_{k_0i}=0$$
  in the tensor product $I_{\cal F}\otimes _{C(I)}I_{\cal F},$ where
  $C(I)$ is the generalized centroid of $I$ and $I_{\cal F}$ is the left
  Martindale ring of quotients of $I.$ It is well-known and it is easy to
  see that $I_{\cal F}=R_{\cal F}$ and $C(I)=C(R).$ Therefore for any
  $x_1\in R_{\cal F}$ we have the identity
  $$\sum_i v_{k_0i}x_1w_{k_0i}=0.$$
  This identity with (21) and (20) implies that the identity
    $$c_{k_0}^{\prime }f(\sum_i a_{i1}x_1a_{i1}^*, u_2(b_2),\ldots
    ,u_n(b_n))=0 \eqno(22)$$
 is valid for each $x_1\in R_{\cal F}.$

Since the values $b_2,\ldots ,b_n$ are arbitrary from $I$, we have an identity
of the form
    $$f(\sum_i a_{i1}x_1a_{i1}^*,u_2(x_2),\ldots ,u_n(x_n))=0, \eqno(23)$$
    where $x_1\in R_{\cal F}, x_2\in I, \ldots ,x_n\in I.$

Now let us fix values 
$x_1\in R_{\cal F}, x_3=b_3\in I, \ldots ,x_n=b_n \in I.$
Then in the same way we obtain
   $$f(\sum_i a_{i1}x_1a_{i1}^*,\sum_i a_{i2}x_2a_{i2}^*,\ldots
    ,u_n(x_n))=0,$$
where $x_1,x_2\in R_{\cal F}, x_3,\ldots x_n \in I.$

Continuing this process we will obtain the following identity on $R_{\cal F}:$ 
    $$f(\sum_i a_{i1}x_1a_{i1}^*,\sum_i a_{i2}x_2a_{i2}^*,\ldots
    ,\sum_i a_{in}x_na_{in}^*)=0, \eqno(24)$$
This is a generalized identity valid in $R_{\cal F}\supseteq R.$ All we need
is to prove that for some $\lambda _1,\ldots, \lambda_n;\rho _1, \ldots
,\rho _n$ this is not a trivial identity. It means that the left hand side
of (24) is not zero in the free product $R_{\cal F} *_CC\langle
x_1,\ldots ,x_n\rangle $
or, in other words, this identity does not follow from the trivial
 generalized identities
$xc=cx,$ where $c\in C.$
Otherwise assume all these identities are trivial.

 Any application
 of a trivial
identity does not change the order of the indeterminates, therefore
all the generalized monomials (i.e. sums of all monomials with fixed order
of sequence of the indeterminates) in the identities (24) should be (trivial)
identities. These generalized monomials have the form
   $$\alpha _{\pi }(\sum_i a_{i\pi (1)}x_{\pi (1)}a^*_{i\pi (1)})
   (\sum_i a_{i\pi (2)}x_{\pi (2)}a^*_{i\pi (2)})\cdots
   (\sum_i a_{i\pi (n)}x_{\pi (n)}a^*_{i\pi (n)}),$$
where $\pi $ is a permutation and
   $$f(x_1,\ldots x_n)=\sum_{\pi }\alpha _{\pi}x_{\pi (1)}\cdots x_{\pi (n)}.$$
Since one of the coefficients $\alpha _{\pi }$ is equal to one
 (let $\alpha _1=1$),
  $$(\sum_i a_{i1}x_1a^*_{i1})(\sum_i a_{i2}x_2a^*_{i2})\cdots
  (\sum_i a_{in}x_na^*_{in})=0 \eqno(25)$$

Let us fix some values of $x_2,\ldots ,x_n$ in $R$ and apply Proposition 4.3
to (25), where we suppose $x=x_1,$ and all coefficients $a_{ki}, \ \ k=2,3,\ldots
p^m$ are zero. We have
   $$(\sum_ia_{i1}\otimes a_{i1}^*)(\sum_i a_{i2}x_2a^*_{i2})\cdots
  (\sum_i a_{in}x_na^*_{in})=0.$$
The set $\{ a_{i1}\} $ is a basis of the ideal $\lambda _1,$
i.e. this is a linearly independent set, therefore
  $$a_{i1}^*(\sum_i a_{i2}x_2a^*_{i2})\cdots
  (\sum_i a_{in}x_na^*_{in})=0$$
for all $a_{i1}^*$ from the dual basis $\{ a_{i1}^*\} $ of the conjugated
ideal $\rho _1.$ This implies that
  $$\rho _1(\sum_i a_{i2}x_2a^*_{i2})\cdots
  (\sum_i a_{in}x_na^*_{in})=0.$$
Since the pair $(\lambda _1,\rho _1)$ was chosen in an arbitrary way,
  $$(\sum_{\rho ^*\simeq\ a\ left\ ideal\ of\ {\cal B}}\rho )
  (\sum_i a_{i2}x_2a^*_{i2})\cdots
  (\sum_i a_{in}x_na^*_{in})=0. \eqno(26)$$
By Property (Q5) of quasi-Frobenius algebras $$1\in {\cal B}=
(\sum_{\rho ^*\simeq\ a\ left\ ideal\ of\ {\cal B}}\rho ) $$ and
therefore
  $$(\sum_i a_{i2}x_2a^*_{i2})\cdots
  (\sum_i a_{in}x_na^*_{in})=0.$$

Now the evident induction works. The theorem is proved.\

The same proof can be applied also for some cases when the inner
part is not quasi-Frobenius but has enough pairs of conjugated
one-sided ideals. Indeed, let us denote by ${\cal B}_r$ the sum
of all right ideals of a finite dimensional algebra ${\cal B}$
conjugated to left ones.

{\bf 7.2. Lemma.} ${\cal B}_r$ {\it is a two-sided ideal of} ${\cal B}.$

{\bf Proof.} Let $\rho $ be a right ideal such that the dual left
module $\rho ^*={\rm Hom}(\rho ,C)$ is isomorphic to a left ideal
$\lambda .$ If $b\in {\cal B}$ then we have an exact sequence of homomorphisms
of right ideals $\rho \rightarrow b\rho \rightarrow 0.$ The conjugated
sequence has the form $\rho ^*\leftarrow (b\rho )^*\leftarrow 0,$
therefore the right ideal $b\rho $ has a conjugated module $(b\rho )^*$
which is isomorphic to a left subideal of $\lambda \simeq \rho ^*.$
Thus $b\rho \subseteq {\cal B}_r$ and ${\cal B}_r$ is a two-sided ideal.
The lemma is proved.

In the same way one can define an ideal ${\cal B}_l$ --- the sum of all
left ideals conjugated to right ones.

{\bf 7.3. Theorem.} {\it Let $L$ be a finite dimensional restricted
differential Lie $C$-algebra of $R$-continuous derivations of a prime
ring $R$ of positive characteristic $p>0.$ If the algebra of
constants $R^L$ satisfies a multilinear polynomial identity of degree
$n$ and ${\cal B}(L)^n_r\neq 0,$ then $R$ is a GPI-ring.}

{\bf Proof.} In the same way as in the proof of Theorem 7.1 we have
identities (24). If all of these identities are trivial then we also have
the identities (26)
 which can be written in the form
   $${\cal B}(L)_r(\sum a_{i2}x_2a_{i2}^*)\cdots (\sum
   a_{in}x_na_{in}^*)=0. \eqno(27)$$

If $b$ is an arbitrary element  of ${\cal B}(L),$ then
$b(\sum_i a_{ik}x_ka_{ik}^*)=(\sum _i a_{ik}x_ka_{ik}^*)b.$
Therefore for $b\in {\cal B}(L)_r,$ identity (27) implies
   $$(\sum a_{i2}x_2a_{i2}^*)\cdots (\sum
   a_{in}x_na_{in}^*)b=0. \eqno(28)$$
By Proposition 4.3 we have
  $$(\sum a_{i2}\otimes a_{i2}^*)\cdots (\sum
   a_{in}x_na_{in}^*)b=0,$$
as in the proof of Theorem 7.1 we have
  $${\cal B}(L)_r(\sum a_{i3}x_3 a_{i3}^*)\cdots (\sum
   a_{in}x_na_{in}^*)b=0,$$
thus
  $$(\sum a_{i3}x_3 a_{i3}^*)\cdots (\sum
   a_{in}x_na_{in}^*){\cal B}(L)_r^2=0.$$
Now the evident induction implies ${\cal B}(L)_r^n=0.$ Hence
 one of the GPI's (24) is not trivial. The theorem is proved.

In a symmetrical way one can prove that the condition ${\cal B}(L)^n_l\neq 0$
also implies that one of the identities (24) is not trivial. It can be proved
that ${\cal B}^n_r=0$ iff ${\cal B}^n_l=0:$

{\bf 7.4. Proposition.}  {\it Let ${\cal B}$ be a finite dimensional algebra.
Then all $n+1$ conditions ${\cal B}_r^k{\cal B}_l^{n-k}=0, \ \ k=0,\ldots ,n$
are equivalent to each other.}

{\bf Proof.} It is enough to show that the conditions for $k$
 and $k+1$ are equivalent. The condition ${\cal B}_r^k{\cal B}_l^{n-k}=0$
 is equivalent to
  ${\cal B}_r^k{\cal B}_l^{n-k-1}\lambda =0$ for
 all pairs of conjugated ideals $\rho , \lambda .$  Since the form
 $(\ ,\ ):\rho \times \lambda \rightarrow C$ is 
 nondegenerate the last condition
 for given $\lambda , \rho $ is equivalent to
 $(\rho ,{\cal B}_r^k{\cal B}_l^{n-k-1}\lambda )=0.$ By
 associativity of the form this is equivalent to
 $(\rho {\cal B}_r^k{\cal B}_l^{n-k-1},\lambda )=0$ and since
 the form is nondegenerate
 this is equivalent to $\rho {\cal B}_r^k{\cal B}_l^{n-k-1}=0.$
 The last conditions for all pairs of conjugated ideals $\lambda , \rho $
 are equivalent to ${\cal B}_r^{k+1}{\cal B}_l^{n-k-1}=0.$ The proposition
 is proved.

Now it is a question of interest whether the condition ${\cal B}(L)^n_r=0$
implies that all identities (24) are trivial generalized polynomial
identities. The answer is yes:

{\bf 7.5. Proposition.} {\it If under the conditions of theorem {\rm 7.3}
 \ ${\cal B}(L)_r^n=0,$ then all identities {\rm (24)} are trivial.}

{\bf Proof.} It is enough to show that all the generalized monomials (25)
are trivial identities. We will prove by inverse induction on $k$
 that for arbitrary
$b_1, \ldots ,b_k\in {\cal B}(L)_r$ the generalized polynomial
  $$(\sum _ia_{i\ k+1}x_{k+1}a_{i\ k+1}^*)\cdots
  (\sum _ia_{in}x_na_{in}^*)b_kb_{k-1}\cdots b_1=0 \eqno(29)$$
is a trivial generalized identity.

If $k=n$ then (29) has the form $b_nb_{n-1}\cdots b_1=0$ that is a trivial
identity as ${\cal B}(L)^n_r=0.$

Assume that (29) is a trivial identity. The identities
$$b(\sum
_i a_{is}xa_{is}^*)=(\sum
_i a_{is}xa_{is}^*)b, \ \ b\in {\cal B}(L) \eqno(30)$$
are trivial generalized polynomial identities (as well as any
linear generalized identity). Let $b_k=a_{ik}^*,$ then from (29) and (30)
we have the following trivial identity
  $$a_{ik}^*(\sum _ia_{i\ k+1}x_{k+1}a_{i\ k+1}^*)\cdots
  (\sum _ia_{in}x_na_{in}^*)b_{k-1}\cdots b_1=0.$$
  Multiplication of this equality on the left by $a_{ik}x_k$ and summation
over $i$ gives the equality (29) with a smaller $k.$ The proposition is 
proved.\

\section{Semiprime PI-rings of constants}

In this secton we will prove under the conditions of Theorem 7.1,
that if
the  ring of constants  $R^L$ is a semiprime PI-ring, then $R$ is also PI.

{\bf 8.1. Theorem.} {\it Let $L$ be a finite dimensional restricted
differential Lie $C$-algebra of $R$-continuous derivations of a prime
ring $R$ of positive characteristic $p>0.$ Suppose that the inner associative
part ${\cal B}(L)$ of $L$ is quasi-Frobenius. If the ring of constants $R^L$
is a semiprime PI-ring, then $R$ is PI.}

{\bf Proof.} By Theorem 7.1 the ring $R$ satisfies a generalized polynomial
identity. Moreover all generalized polynomial identities (24) hold
in its left Martindale ring of quotients $R_{\cal F}.$ In particular they hold
in the central closure $RC\subseteq R_{\cal F}$ of the ring $R.$
By the Martindale structure theorem
 [Ma69] this central closure has an idempotent
$e,$ such that $D=eRCe$ is a skew field of  finite dimension
over $C.$ (Note that formally  Martindale
theorem can be applied only if the
coefficients of the identity belong to $R.$ In our case they belong to
$R_{\cal F}$ but may not belong to $R.$ Nevetheless Martindale's original
proof is correct for our case too; see, for instance, [Kh91, Theorem 1.13.4]
or the special investigation in [La86].)

Thus, by the Martindale theorem, $RC$ is a primitive ring with a nonzero socle.
The N. Jacobson structure theorem [Ja64] shows that $RC$ is a dense subring in
the finite topology in the complete ring ${\cal E}$ of linear transformations
of the left space $V=eRCe$ over the skew field $D.$

Moreover, the left Martindale quotient ring $(RC)_{\cal F}$ is equal to
${\cal E}$ (see, [Ha82, Lemma 1.1] and [Ha87, Remark 4.9] or
[Kh91, Theorem 1.15.1]). It is easy to see
that $R_{\cal F} \subseteq (RC)_{\cal F}={\cal E}.$ (Indeed, if $q\in R_{\cal
F}$ and $Iq\subseteq R$ for a nonzero ideal $I$ of $R,$ then we can extend
$q$ to the ideal $IC$ of $RC$ by the obvious formula
$(\sum i_{\alpha }c_{\alpha })q=\sum (i_{\alpha }q)c_{\alpha }.$
This is well-defined. Indeed, if $\sum i_{\alpha }c_{\alpha }=0$ 
and $J$ is a nonzero
 ideal of $R$
such that $Jc_{\alpha }\subseteq R$ then $\sum (jc_{\alpha })i_{\alpha }=0$
for all $j\in J.$ Therefore $\sum (jc_{\alpha })(i_{\alpha }q)=0;$ i.e.
$J(\sum c_{\alpha }(i_{\alpha }q))=0$ and $\sum (i_{\alpha }q)c_{\alpha }=0.$)
Now all the coefficients of (24) belong to
${\cal E}$ and since addition and multiplication are continuous in the finite
topology, the identities (24) hold in ${\cal E}.$ (Here one can use
also Corollary 2.3.2 from [Kh91] which
allows us to extend identities from $RC$ to
$(RC)_{\cal F}.$)

Now we are going to prove that the space $V$ is finite dimensional over $D.$
In that case the dimension of ${\cal E}$ over $C$ will also be finite:
$d=\dim _C{\cal E}= (\dim _DV)^2 \cdot \dim _CD$ and ${\cal E}$ (and 
therefore $R$),
like any $d$-dimensional algebra, will
satisfy the standard polynomial identity:
$$S_d(x_1, \ldots , x_{d+1})\equiv \sum (-1)^{\pi }x_{\pi (1)}\cdots x_{\pi
(d+1)}=0.$$

On the contrary, suppose that $V$ has infinite dimension $\dim V=\beta .$
Let $M$ be the set of all linear transformations whose rank is less then $\beta .$
(Recall that the {\it rank } of a transformation $l$ is the
 dimension over $D$ of
its image.) It is well-known that $M$ is a maximal ideal of ${\cal E}.$
So the factor ring $\bar{\cal E}={\cal E}/M$ is a simple ring with a unit.

{\bf 8.2. Lemma.}  { \it The ring $\bar{\cal E}$ is not Artinian.}

{\bf Proof.} Let $\{ e_i, i\in I\} $ be a basis of $V$ over $D.$ and
$$I_1\supset I_2\supset \ldots \supset I_n\supset \ldots $$ be a chain of
subsets such that $|I_k \setminus I_{k+1}|=\beta ,$ and let
$$A_n=\{ l\in {\cal E}:e_il=0 \ \  \forall i \in I\setminus I_n\} .$$
Then
$$ (A_1+M)/M\supset A_2+M/M\supset \ldots \supset A_n+M/M\supset \ldots $$
is an infinite descending chain of right ideals of $\bar{\cal E}.$

Indeed, if $A_n+M=A_{n+1}+M,$ then for the transformation $w,$ defined by
\[ e_iw=\left\{ \begin{array}{ll}
              e_i & \mbox{if $i\in I_n\setminus I_{n+1}$} \\
              0  & \mbox{otherwise}
              \end{array}
              \right. , \]
we should get a presentation $w=a+m,$ where $a\in A_{n+1}, \ m\in M.$
Let $V_1$ be a subspace generated by $\{e_i:i\in I_n\setminus I_{n+1} \} .$
Then $V_1=V_1w\subseteq V_1a+V_1m=V_1m.$ However, $\dim _DV_1=\beta ,$
while $\dim _DV_1m\leq \dim Vm<\beta ,$ which is a contradiction. The lemma is
proved.\

{\bf 8.3. Lemma.} {\it The ring $\bar{\cal E}$ does not satisfy a non trivial
generalized polynomial identity.}

{\bf Proof.} Like any simple ring with a unit, the ring $\bar{\cal E}$ is
primitive. If it satisfies a GPI, then by the S.A. Amitsur structure theorem
[Am65] it has a nonzero socle $S,$ which is a two-sided ideal and
therefore $S=\bar{\cal E}.$ In N. Jacobson presentation of $\bar{\cal E}$
as a dense ring of linear transformations, the socle consists of all
transformations of finite rank. This means that the unit has finite
rank and therefore the space has finite dimension. Thus $\bar{\cal E}$
is the ring of all linear transformations of a finite dimensional space
over a skew field. In particular $\bar{\cal E}$ is Artinian; this is a
 contradiction
to Lemma 8.2. The lemma is proved.

Let us consider now identities (24). We have seen that all these identities
hold in ${\cal E}.$ If we apply the natural homomorphism
$\varphi :{\cal E}\rightarrow \bar{\cal E}={\cal E}/M$ we obtain
the following identities of the ring $\bar{\cal E}$
$$f(\sum _i\bar{a}_{i1}x_1\bar{a}^*_{i1},\ldots ,\sum
_i\bar{a}_{in}x_n\bar{a}^*_{in})=0, \eqno(31)$$
where $\bar{a}=\varphi (a)=a+M.$

By Lemma 8.3  all we need is to prove that
 one of the identities (31) is a nontrivial
GPI of $\bar{\cal E}.$

First of all we have to calculate the generalized centroid of $\bar{\cal E}.$
As $\bar{\cal E}$ is a simple ring with a unit, it equals its left Martindale
quotient ring and therefore the generalized centroid is equal to the center.

{\bf 8.4. Lemma.}  {\it The center of
 $\bar{\cal E}$ is canonically isomorphic
to $C, \ \ C(\bar{\cal E})=\varphi (C).$}

{\bf Proof.} See [Ro58, Corollary 3.3]. 
We will need the following result which gives a criterium
for determining when the
 ring of constants is semiprime (see Theorem 5.1 [Kh95]).

{\bf 8.5. Theorem.} {\it Under the conditions of theorem 8.1, the ring of
constants is semiprime if and only if ${\cal B}(L)$ is differentially
semisimple, i.e. it has no nonzero differential (with respect to
action of $L$) ideals with zero
multiplication or, equivalently, it is a sum of a finite number of
differentially simple algebras.}

By this theorem we have that in our situation the algebra ${\cal B}(L)$
is differentially semisimple.

{\bf 8.6. Lemma.} {\it The ideal $M$ is a differential ideal with
respect to $L,$ i.e. $M^{\mu }\subseteq M$ for each $\mu \in L.$}

{\bf Proof.} Note that $M$ is a differential ideal with respect to
each derivarion of ${\cal E}.$ Indeed, if $l\in M$ than $l$ is a
transformation of rank less then $\beta $ and the projection
$e:V\rightarrow {\rm im} l$ also has rank less than $\beta ,$
in which case $l=le.$ We have $l^{\mu }=l^{\mu }e+le^{\mu }\in M$
for each derivation $\mu \in {\rm Der}({\cal E}).$

By proposition 1.8.1 [Kh91] any $R$-continuous derivation has a unique
extension to $R_{\cal F}.$ In particular each derivation from $L$
is defined on $RC.$ Again by the same proposition we have that
the elements of $L$ have a unique extensions 
to $(RC)_{\cal F}={\cal E}.$ Thus we have
obtained that the ideal $M$ is differential with respect to $L.$
The lemma is proved.

As a consequence we have that the intersection $M_0=M\cap {\cal B}(L)$ is
a differential
ideal of ${\cal B}(L),$ which is not equal to ${\cal B}(L)$
(it does not contain 1). The left
annihilator $l(M_0)$ of $M_0$ in ${\cal B}(L)$ is also a differential
ideal, therefore $l(M_0)\cap M_0 $ is a differential ideal with zero
mulitiplication. By theorem 8.5, $l(M_0)\cap M_0=0.$ In the same way
the left annihilator of the sum $l(M_0)+M_0$ is zero (it is contained in
$l(M_0)$ and, therefore, has a zero multiplication). Now property (Q1)
of quasi-Frobenius algebras implies that
$l(M_0)+M_0=r(l(l(M_0)+M_0))=r(0)={\cal B}(L)$ and, finally
$${\cal B}(L)=l(M_0)\oplus M_0=e{\cal B}(L) \oplus (1-e){\cal B}(L), \eqno(37)$$
where $e$ is a central idempotent defined by the corresponding decomposition
of the unit $1=e\oplus (1-e).$

Let us return to identities (31). Suppose that in these identities
$\{ a_{ij}\} $ and $\{ a_{ij}^* \} $ are bases of conjugated ideals
$\lambda _j, \rho _j$ contained in $l(M_0).$ In that case the sets
$A_j=\{ \bar{a}_{ij}, i=1, \ldots m \} $ are linearly independent over
the center of $\bar{\cal E}$ (see lemma 8.4). Moreover, the $C$-space
  generated by all possible $a_{ij}^*$'s contains the
unit $e$ of $l(M_0)$ because for each 
conjugated pair of ideals $\lambda ,\rho $
the one-sided ideals $e\lambda ,e\rho $ are also conjugated with 
respect to the same
form (note that $e$ is a central idempotent of ${\cal B}(L)$).
This implies that the linear space over the center of $\bar{\cal E}$
generated by all $\bar{a}_{ij}^*$'s contains the unit $\bar{e}$ of
$\bar{\cal E}.$  This fact allows us to prove that one of the identities
(31) is nontrivial in the same manner as it was done in the end of
the proof of  Theorem 7.1. By Lemma 8.3, Theorem 8.1 is proved.

In this proof we used the fact that the 
inner part ${\cal B}(L)$ is differentially
semisimple and that it has enough pairs of conjugated ideals.
Therefore in the way analogous to Theorem 7.3 we can formulate a
slightly more general result.

{\bf 8.7. Theorem.} {\it Let $L$ be a finite dimensional restricted
differential Lie $C$-algebra of $R$-continuous derivations of a prime
ring $R$ of positive characteristic $p>0.$ Suppose that the inner
part ${\cal B}(L)$ is a direct sum of differentially simple ideals
$${\cal B}(L)=B_1\oplus B_2\oplus \ldots \oplus B_m.$$
If the algebra of constants $R^L$ satisfies a multilinear polynomial
identity of degree $n$ and $(B_i)^n_r\neq 0, i=1,\ldots ,m,$ then
$R$ is a PI-ring.}

The only place where we have used that ${\cal B}(L)$ is quasi-Frobenius
is  decomposition (37). Therefore it is enough to show that
 each differential ideal of
 the direct  sum of differentially simple algebras with units
  is a direct summand.
 If $B=B_1\oplus B_2\oplus \ldots \oplus B_m$ is a direct sum of
 differentially simple algebras then for any differential ideal $A$
 we have that $A B_i$  is a differential ideal of $B_i.$ This implies
 that either $B_i\subseteq A$ or $A B_i=0.$ In the same way
 either $B_i\subseteq A$ or $B_iA=0.$
 Let $l(A)$ be the left annihilator of $A,$ then $l(A)\cap A$ is
 a differential ideal with zero multiplication, so its product
 with each $B_i$ is zero. This is possible only if the intersection
 is zero. In the same way the left annihilator of the sum $l(A)+A$
 has a zero multiplication and therefore it is equal to zero. It means that
 $l(A)+A$ contains all the components $B_i$ and $l(A)\oplus A=B.$
 
\

ACKNOWLEDGMENT

The authors are grateful to Professor Dalit Baum for her help.

\

{\bf REFERENCES}

\

[Am65]. S.A.Amitsur, {\it Generalized polynomial identities and pivotal
monomials}, Trans. Amer. Math. Soc., v.114(1965), 210--216.

[Ba27]. R.Baer, {\it Algebraiche theorie der differentierbaren funktionen
koper}, I. -- Sitzungsberichte. Heidelberger Academia, 1927, 15--32.

[Be64]. A.J.Berkson, {\it The u-algebra of a restricted Lie algebra is
Frobenius}, PAMS, v.15(1964), 14--15.

[CR62]. C.W.Curtis and I.Reiner, {\it Representation Theory of Finite
Groups and Associative Algebras}, Interscience Publishers, New York -- London,
1962.

[Ha82]. M.Hacque, {\it Anneaux fidelement repr\'esent\'es sur leur socle 
droit}, Communications in Algebra, v.10, no.10(1982), 1027--1072.

[Ha87]. M.Hacque, {\it Th\'eorie de Galois des anneaux presque-simples},
Journal of Algebra, v.108, no.2(1987), 534--577.

[Ja64]. N.Jacobson, {\it Structure of Rings}, Amer. Math. Soc. Colloquium.
Publ., Providence, 1964.

[Kh78]. V.K.Kharchenko, {\it Differential identities of prime rings},
Algebra i logika, v.17, no.2(1978), 220--238.

[Kh81]. V.K.Kharchenko, {\it On centralizers of finite dimensional algebras},
Algebra i logika, v.20, no.2(1981), 231--247.

[Kh82]. V.K.Kharchenko, {\it Constants of derivations of prime rings},
Math. USSR Izvestija, v.18, no.2(1982), 381--401.

[Kh91]. V.K.Kharchenko, {\it Automorphisms and Derivations of Associative
Rings}, Kluwer Academic Publishers, v.69(1991).

[Kh95]. V.K.Kharchenko, {\it On derivations of prime rings of 
positive characteristic}, Algebra i logika, 1995, to appear.

[Ko91]. A.N.Korjukin, { \it To a question of bicentralizers in prime rings},
Sib. Mat. Journal, v.32, no.6(1991), 81--86.

[La86]. C.Lanski, {\it A note on GPIs and their coefficients}, Proc. Amer.
Math. Soc. v.98(1986), 17--19.

[Lv93]. I.V.Lvov, {\it  On centralizers of finite dimensional subulgebras
in the algebra of linear transformations}, Third International Algebraic 
Conference, Krasnojarsk, 1993, 213--214.

[Ma69]. W.S.Martindale, {\it Prime rings satisfying a generalized polynomial
identity}, Journal of algebra, v.12, no.4(1969), 576--584.

[Pa87]. D.S.Passman, {\it Prime ideals in enveloping rings}, Trans. Amer.
Math. Soc., v.302(1987), 535--560.

[Pi86]. Piers Dos Santos, {\it Derivationes des anneaux semi-premiers I},
Comm. in algebra, v.14, no.8(1986), 1523--1559.

[Po83]. A.Z.Popov, {\it On derivations of prime rings}, Algebra i Logika,
v.22, no.1(1983), 79--92.

[Ro58]. A.Rosenberg, {\it The structure of the infinite general linear 
group}, Ann. Math. ser.2, 68(1958), 278--294
\end{document}